\documentclass{ifacconf}

\usepackage{graphicx}      
\usepackage{natbib}        
\usepackage{amsmath}
\usepackage{amssymb}
\usepackage{latexsym}
\usepackage{mathrsfs}
\usepackage{amscd}
\usepackage{sidecap}
\usepackage{graphicx}
\usepackage{color}
\usepackage{stmaryrd}
\usepackage{algorithmic}
\usepackage{algorithm}
\graphicspath{{pics/}}
\def\bigK{ {\mathcal {K}}}

\def\U{ {\mathcal {U}}}

\def\QF{\mathcal{P}}

\def\mrho{\mu}
\newcommand{\R}{\mathbb R}
\newcommand{\J}{\mathcal J}
\newcommand{\N}{\mathbb N}
\newcommand{\dt}{\Delta t}

\newcommand{\epsi}{{\varepsilon}}

\def\be#1\ee{\begin{equation}#1\end{equation}}

\newcommand{\bx}{\mathbf{x}}
\newcommand{\bu}{\mathbf{u}}

\renewcommand{\theequation}{\arabic{section}.\arabic{equation}}
\renewcommand{\thetable}{\arabic{section}.\arabic{table}}
\renewcommand{\thefigure}{\arabic{section}.\arabic{figure}}
\newcommand{\bq}{\begin{equation}}
\newcommand{\eq}{\end{equation}}

\newtheorem{remark}{Remark}[section] 

\def\bqa{\begin{eqnarray}}
\def\eqa{\end{eqnarray}}
\def\dt{\Delta t}

\def\S{\sigma}

\newcommand{\bd}{\begin{displaymath}}
\newcommand{\ed}{\end{displaymath}}
\newcommand{\ba}{\begin{eqnarray}}
\newcommand{\ea}{\end{eqnarray}}

\def\L{\mathcal L}
\def\R{\mathbb{R}}
\def\N{\mathbb{N}}
\def\R{\mathbb{R}}

\def\pa{\partial}

\def\epsi{ \varepsilon}

\newcommand{\lt}{\left}
\newcommand{\rt}{\right}

\newcommand{\K}{\mathcal{K}}

\usepackage{graphicx}      
\usepackage{natbib}        
\usepackage{ams math}
\begin{document}
\begin{frontmatter}

\title{A Boltzmann approach to mean-field sparse feedback control \thanksref{footnoteinfo}} 

\thanks[footnoteinfo]{Giacomo Albi and Massimo Fornasier acknowledge the support of the ERC-Starting Grant HDSPCONTR "High-Dimensional Sparse Optimal Control". Dante Kalise acknowledges the support of the ERC-Advanced Grant OCLOC "From Open-Loop to Closed-Loop Optimal Control of PDEs".}

\author[First]{Giacomo Albi} 
\author[Second]{Massimo Fornasier} 
\author[Third]{Dante Kalise}

\address[First]{Department of Mathematics, TU M\"unchen, Boltzmannstr. 3, Garching bei M\"unchen,  D-85748, Germany (e-mail: giacomo.albi@ma.tum.de).}
\address[Second]{Department of Mathematics, TU M\"unchen, Boltzmannstr. 3, Garching bei M\"unchen,  D-85748, Germany (e-mail: massimo.fornasier@ma.tum.de)}
\address[Third]{Johann Radon Insitute for Computational and Applied Mathematics, Austrian Academy of Sciences, Altenbergerstr. 69, A-4040 Linz, Austria, (e-mail: dante.kalise@oeaw.ac.at)}

\begin{abstract}                
We study the synthesis of optimal control policies for large-scale multi-agent systems. The optimal control design induces a parsimonious control intervention by means of $\ell_1$, sparsity-promoting control penalizations. We study instantaneous and infinite horizon sparse optimal feedback controllers. In order to circumvent the dimensionality issues associated to the control of large-scale agent-based models, we follow a Boltzmann approach. We generate (sub)optimal controls signals for the kinetic limit of the multi-agent dynamics, by  sampling of the optimal solution of the associated two-agent dynamics. Numerical experiments assess the performance of the proposed sparse design.
\end{abstract}


\end{frontmatter}
\section{Introduction}
Multi-agent dynamical systems (MAS) naturally arise in the mathematical modeling of social dynamics in a wide spectrum of applications: animal behavior, cellular aggregation, opinion dynamics, and human crowd motion, among many others (\cite{CDFSTB03,cristiani2014book}). So far, it has been of utmost interest to study different collective behavior phenomena such as clustering or consensus emergence without external forces. Depending on the degree of cohesiveness of the initial configuration of the agents and the strength of their interaction, dynamical patterns may arise naturally by self-organization. However, if self-organization is not sufficient to enforce a stable pattern, collective behavior can be induced by means of exogenous interventions.
In this paper, we study the design of control actions which are able to steer a MAS towards a prescribed consensus regime. We address this challenge by means of optimal control techniques, thus minimizing an energy measure of both the control and the state of the system, constrained to the multi-agent dynamics. MAS are naturally represented as large-scale systems of $N$ coupled nonlinear differential equations of the form
\begin{align*}
 \dot{x}_i(t)= \frac{1}{N}\sum_{j=1}^N P(x_i,x_j)(x_j-x_i)+u_i(t), \qquad i=1,\ldots,N,
\end{align*}
where $x_i(t)$ represents the state of the $i-$th agent, and the binary kernel $P(\cdot,\cdot)$ encompasses interaction rules between agents such as attraction, repulsion, or alignment. The term $u_i(t)$ corresponds to a dynamical external action, which we design in an optimization-based framework.

As the number of agents increases, the complexity of casting an optimal control problem becomes prohibitively expensive, a phenomenon often referred as Bellman's curse of dimensionality (\cite{BELL}). In order to circumvent this difficulty, we follow a multiscale approach. By borrowing a leaf from statistical mechanics, the mean field approximation of a multi-agent system replaces the microscopic representation of the state by an agent density function, which evolves according to a nonlinear, nonlocal kinetic transport equation of the form
\begin{align}
\pa_t \mrho + \nabla \cdot \lt(\lt(\QF[\mrho] + u\rt)\mrho \rt) = 0,
\end{align}
where $\mu=\mu(x,t)$ is the density of agents at time $t$ at state $x$, and $u=u(x,t)$ is the mean-field realization of the external forcing. The interaction force $\mathcal P$ is given by
\begin{equation}\label{eq:convolution}
\QF[\mrho](x) = \int_{\R^d} P(x,y)(y-x)\mrho(y,t)\,dy\,.
\end{equation}
With an adequate characterization of such a  controller, the optimal design problem can be understood as a fluid flow control problem. Furthermore, the so-called Boltzmann approach yields an approximation of the mean-field dynamics by means of an iterative sampling of 2-agent microscopic dynamics (binary dynamics). This same principle allows us to generate control signals for the mean-field model by means of solving optimal control problems associated to the binary dynamics.

The topic of emergent collective behavior in MAS has been linked the study of pattern formation and self-organization phenomena (\cite{2007TAC,CFRT}), and to recent developments covered within the area of Mean Field Games (\cite{huang2006,LL}). This latter approach includes an optimal and decentralized decision process, as for instance in the financial market, and the emphasis is in the characterization of Nash equilibria. We follow a different approach, enforcing consensus by optimizing the intervention of a centralized external policy maker endowed with limited resources. This approach has been already studied for microscopic dynamics in \cite{CFPT,BFK,borzi2015M3ASa}, and at the mean-field level in \cite{BFY,FS13,BFRS}, among others. In \cite{ABCK15}, we have developed an analytical and computational multiscale optimal control approach for the control of mean-field dynamics through the inclusion microscopic leaders, with applications in crowd motion evacuation. More recently, in \cite{ACFK16} we introduced a mean-field control hierarchy, where optimal feedback controllers are computed for a binary system of particles, and its action is inserted in the mean-field dynamics, regulating the density evolution towards a target. In this work, based on this modeling/control approach, we focus on the design features of the optimal design, more specifically, in the formulation of sparse optimal control problems, where parsimonious control action is enforced by means of a non-smooth $\ell_1$-norm penalization of the control term. 

Sparse control of microscopic and mean-field MAS has been previously studied in the context of necessary optimality conditions (\cite{CFPT,MR3268059}). In this paper we differ from this approach by addressing the design of sparse optimal feedback controllers based on a dynamic programming formulation. We study two particular limit cases: the design of instantaneous controllers, and infinite horizon optimal control. Instantaneous optimal controllers for mean-field MAS have been studied in \cite{AHP,APZa} for $\ell_2$-norm control penalization, and their appeal resides in the fact that the arbitrary small time frame in which the control is synthesized allows the computation of closed-form feedback solutions, which are real-time implementable. However, an instantaneous controller can guarantee asymptotic stability of the closed-loop dynamics only in very specific cases, and will require continuous control action unless sparsity in time is enforced. Instead, the infinite horizon optimal feedback control depends on the solution of a Hamilton-Jacobi-Bellman equation, producing a nonlinear feedback mapping which is suitable for real-time control and can guarantee asymptotic stability. Infinite horizon sparse optimal feedback control has been recently addressed in (\cite{KKK16,KKR16}) for general nonlinear dynamics.

This paper is organized as follows. In Section 2 we state the optimal control formulation for multi-multi-agent dynamics. In Section 3, we present two different sparse optimal control problems for 2-agent models. Section 4 is devoted to the Boltzmann approach which allows us to compute (sub)optimal controllers for the mean-field approximation of the MAS by means of iterative sampling of the controlled binary model. Finally, Section 5 presents numerical experiments related to optimal control of opinion dynamics which illustrate the different features of the proposed designs.

\section{Microscopic optimal control}
We consider a population of $N$ agents represented by $x_i(t)\in\R^d$, evolving according to
\begin{align}
 \frac{dx_i}{dt}&= \frac{1}{N}\sum_{j=1}^N P(x_i,x_j)(x_j-x_i)+u_i(t)\,,\label{eq:MAS1}\\
 x_i(0)&=x_0\,,\qquad i=1,\ldots,N\,,\label{eq:MAS2}
\end{align}
where $P(\cdot,\cdot): \R^d\times\R^d\longrightarrow\R^d$ is a Lipschitz-continuous communication function, and the control variables $u_i(t)\in\U=\{u(t): \R_+\longrightarrow U\}$, with $U$ a compact subset of $\R^d$. We denote by $\bx(t)=(x_1(t),\ldots,x_N(t))^t$, and $\bu(t)=(u_1(t),\ldots,u_N(t))^t$.
The controllers are obtained as the solution of the following optimal control problem
\begin{align}\label{eq:cost}
\underset{\bu(\cdot)\in\U^N}{\min} \J(\bu(\cdot);\bx_0):=\int_0^T e^{-\lambda t}\ell(\bx(t),\bu(t))\,dt\,,
\end{align}
with a positive discount factor $\lambda>0$, subject to system dynamics \eqref{eq:MAS1}-\eqref{eq:MAS2}. The running cost $l(\bx,\bu)$ is of the form
\begin{align}
\ell(\bx,\bu):=\frac{1}{N}\|\hat\bx-\bx\|_2^2+\gamma\|\bu\|_1,
\end{align}
with $\gamma>0$, $\hat\bx\in\R^{d\times N}$ a desired reference state, and
\begin{align}
\|\bx\|_p=\sum_{i=1}^N|x_i|^{1/p}\,
\end{align}
where $|\cdot|$ stands for the $d$-dimensional Euclidean norm. For the sake of simplicity, in the following we restrict our analysis to the case $d=1$, although the presented methodology is directly applicable to multidimensional agent systems. We enforce sparsity in the design by introducing an $\|\cdot\|_1$ penalization term for the control. 

We shall  assume that system dynamics have been discretized in time with a first-order approximation
\begin{align}
x_i^{k+1}&=x_i^k+\dt\left(\frac{1}{N}\sum_{j=1}^N P(x_i^k,x_j^k)(x_j^k-x_i^k)+u_i^k\right)\,,\label{eq:MASdt}
\end{align}
where $x^k=x(k\dt)$, with $k\in\N$ and a time discretization parameter $\dt>0$. The cost functional \eqref{eq:cost} is discretized accordingly
\begin{align}
\J_{\dt}(\bu;\bx^0):=\sum_{k=0}^{N_T} \beta^k\ell(\bx^k,\bu^k)\,,
\end{align}
with $N_T\dt=T$, and $\beta=e^{-\lambda\dt}$.
For $\dt>0$, we will study the instantaneous control proble ($N_T=1$), and the infinite horizon control problem ( $N_T\to\infty$). In both cases, we shall focus on solutions of \eqref{eq:cost} which can be expressed in feedback form
\begin{align}
\bu^*(t)=\K(\bx(t))\,,
\end{align}
i.e. controllers which can be computed solely based on the information of the current state of system. This global approach characterizes the optimal controller in terms of the (possibly) nonlinear feedback mapping $\K$, which is computed through dynamic programming. Since the dynamic programming approach is limited to low-dimensional dynamics, in the following section we study the control problems for two-particle systems, which will latter generate a sparse feedback controller for the large-scale MAS.

\section{Binary sparse control}
We focus our analysis on the optimal control problem when $N=2$. In this case, denoting by $\bu_{ij} = (u_i,u_j)$ and $\bx_{ij}  = (x_i,x_j)$, we have the following binary control problem:
\begin{align}\label{bincost}
\underset{\bu_{ij}\in\U^2}{\min}\;\;  \J_{\dt}(\bu_{ij} ;\bx_{ij}^0):=\sum_{k=0}^{N_T} \beta^k\ell(\bx_{ij} ^k,\bu_{ij} ^k)\,
\end{align}
subject to the two-agent model
\begin{equation}\begin{split}\label{binary}
x_i^{k+1}&= x_i^k + \frac{\dt}{2} P^k_{ij}(x^k_j-x^k_i) + \Delta t\, u^k_i\,,\\
x_j^{k+1}&= x_j^k + \frac{\dt}{2} P^k_{ji}(x^k_i-x^k_j) + \Delta t\, u^k_j\,,
\end{split}\end{equation}
where $P^k_{ij}:=P(x^k_i,x^k_j)$.

\subsection{Instantaneous sparse control}
The instantaneous control corresponds to the shortest nontrivial prediction horizon, i.e. $N_T=1$. In this case, the control problem is further simplified to
\begin{align}
\underset{\bu_{ij}\in U^2}{\min}\;\;\frac{\beta}{2}\|\hat\bx-\bx_{ij}^{1}\|^2+ \gamma\|\bu_{ij}\|_1 \,,
\end{align}
with $\bx_{ij}^1$ given by eq. \eqref{binary}. Due to the linear dependence of $\bx_{ij}^1$ with respect to $\bu_{ij}$, the minimizer $\bu_{ij}^*$ corresponds to the soft-thresholding operator ({\cite{DJ95}), and is given componentwise by
\begin{align}\label{IC}
u^*_{i}=\Pi_{U}(\mathbb{S}^1_{\gamma,\lambda}(\xi_{i}))\,,
\end{align}
where $\xi_{i}=(\hat{x}_i-x_i^k - \frac{\dt}{2}P^k_{ij}(x^k_j-x^k_i))/\dt$, the operator $\mathbb{S}^1_{\gamma,\lambda}(\xi)$ is defined as 
\[
\mathbb{S}^1_{\gamma,\lambda}(\xi):= \begin{cases}
\left(1-\frac{\bar\gamma}{|\xi|}\right)\xi,&\qquad |\xi|>\bar\gamma\,, \\
0&\qquad otherwise
\end{cases}
\]
with $\bar\gamma=\gamma(\beta \dt^2)^{-1}$, and $\Pi_{U}$ is the projection onto $U$. The expression for $u_{j}^*$ follows analogously. This procedure generates an optimal control in feedback form, i.e. at a given discrete instant $k$, the instantaneous optimal action is a nonlinear mapping only depending on the current state $\bx^k$ and model parameters.
\subsection{Infinite horizon sparse control}\label{sec:HJB}
A more complex feedback synthesis can be performed by considering an infinite prediction horizon, i.e., $N_T=\infty$. In this case, the optimal feedback controller is obtained through dynamic programming. If we define the value function associated to the infinite horizon discrete cost \eqref{bincost} as
\begin{align}\label{eq:Val}
V(\bx_{ij}^0) := \underset{\bu_{ij}\in\U^2}{\inf}\sum_{k=0}^{\infty}  \beta^k \ell(\bx_{ij}^k,\bu_{ij}^k)\,,
\end{align}
then it is well-known that the application of the Dynamic Programming Principle (\cite{BELL}) with the discrete time dynamics \eqref{binary} characterizes the value function as the solution of the Bellman equation
\begin{equation}
\begin{aligned}\label{eq:DP}
V(\bx_{ij}) & = \underset{\bu_{ij}\in\U^2}{\min}\left\{\beta V(\bx_{ij}^+(\bu_{ij}))+\dt\ell(\bx_{ij},\bu_{ij})\right\}\,,
\end{aligned}
\end{equation}
where $\bx_{ij}^+(\bu_{ij})$ denotes a one-step update of \eqref{binary} departing from $\bx_{ij}$ with control action $\bu_{ij}$.
Once this functional equation has been solved for $V$, the optimal feedback controller is given by the nonlinear mapping
\begin{equation}\label{eq:ochj}
\bu^*_{ij}=\underset{\bu_{ij}\in U^2}{\arg\min}\left\{\beta V(\bx_{ij}^+(\bu_{ij}))+\dt\ell(\bx_{ij},\bu_{ij})   \right\}\,.
\end{equation}
In the following, we briefly review the so-called Boltzmann approach for our setting in order to, upon feedback controllers for the binary dynamics, generate control actions for the large-scale MAS.

\section{A Boltzmann approach to mean-field sparse control}
For a large ensemble of agents, the microscopic optimal control problem  \eqref{eq:MAS1}--\eqref{eq:cost} is well-approximated by the following mean-field optimal control problem,
\begin{align}\label{eq:MFJ}
&\min_{u\in\U_{\infty}}\int\limits_{0}^{T}\int\limits_{\mathbb{R}^d}e^{-\lambda t} \L(x,t,\mu(x,t),u(x,t))d\mu(x,t) \ dt,
\end{align}
constrained to the mean-field MAS,
\begin{align}\label{eq:MFC}
\pa_t \mrho + \nabla \cdot \lt(\lt(\QF[\mrho] + u\rt)\mrho \rt) = 0,
\end{align}
where $\mu=\mu(x,t)$ represent the agents' density distribution, evolving from the initial data $\mu(x,0) = \mu^0(x)$. The operator $\mathcal{P}[\mu](x,t)$ is defined as in \eqref{eq:convolution} and $u=u(x,t)$ is the optimal control in a suitable space $\U_{\infty}$. The mean-field running cost $\L$ is defined accordingly to the finite dimensional cost $\ell$. For further details on this type of problems we refer to \cite{FS13}. In general the solution of the variational problem  \eqref{eq:MFJ}--\eqref{eq:MFC}  is challenging problem and computationally demanding due to the nonlinear and nonlocal character of the system dynamics.  We propose an alternative solution procedure, introducing a Boltzmann-type equation  to model the evolution of a system of agents ruled by a binary interactions. This type of description will furnish a suboptimal solution to the mean-field optimal control problem at reduced cost. We briefly review this approach, which is thoroughly developed in \cite{AHP,ACFK16}.

For $\mu=\mu(x,t)$ denoting the kinetic probability density of agents in position $x\in\R^d$ at time $t\geq0$, the time evolution of the density $\mu$ is given as a balance between bilinear gains and losses of the agents' position, due to the following constrained binary interaction,
\begin{equation}
\begin{aligned}\label{eq:bin}
x^*= \,& x        + \alpha P(x,y)(y-x)   + \alpha S_\alpha(x,y),   
\\
y^* =\, & y+ \alpha P(y,x)(x-y)         +  \alpha S_\alpha(y,x),
\end{aligned}
\end{equation}
where $(x^*,y^*)$ are the post-interaction states and the parameter $\alpha$ measures the strength of the interactions, given by  $P(x,y)(y-x)$, and the feedback $S_\alpha(x,y)$ indicating the forcing term due to the control dynamics.
We remark that such dynamics is equivalent to the expression in \eqref{binary} for  $\alpha = \Delta t/2$,$ S_\alpha(x_i,x_j) =  2 u^*_{i}(x_i,x_j)$ and equivalently for  $ S_\alpha(x_j,x_i)$.

We now consider a kinetic model for the  evolution of  the density $\mu$ ruled by the following Boltzmann-type equation 
\begin{align}\label{eq:Boltz}
\pa_t \mu(x,t) = Q_{\alpha}(\mu,\mu)(x,t),
\end{align}
with the interaction operator $Q_{\alpha}(\mu,\mu)$ accounting for the loss and gain of agents in position $x$ via
\begin{align}\label{eq:QBoltzB}
Q_{\alpha}(\mu,\mu)(x) = \eta\int_{\R^d}\left(\frac{1}{\mathbb{J}_\alpha}\mu(x_*) \mu(y_*) - \mu(x)\mu(y)\right)dy,
\end{align}
where $(x_*,y_*)$ are the pre-interaction positions that generate positions $(x,y)$ after the interaction. 
The bilinear operator $Q_\alpha(\cdot,\cdot)$ includes $\mathbb{J}_\alpha$, representing the Jacobian of the transformation  $(x,y)\to(x^*,y^*)$, described by \eqref{eq:bin}, and $\eta>0$ represents a constant interaction rate. We refer to  \cite{PTa,T}, for further generalization and discussion on this type of models.

We want to derive a more regular operator out of the Boltzmann-type of interaction. To this end, we introduce the so called {\em quasi-invariant limit},  whose basic mechanism is to consider a regime where interactions strength is low and frequency is high (\cite{VILL,CPT,T}). In this setting, the following result holds. 
\begin{thm}\label{thm:grazing}
For $\alpha\geq0$ and $t\geq0$, we assume $S_\alpha(\cdot,\cdot)$, $P(\cdot,\cdot)\in L^2_{loc}$ and for $\alpha\to 0$ we assume $S_\alpha(x,y)\to K(x,y)$. Then we consider a weak solution $\mu$ of equation \eqref{eq:Boltz} with initial datum $\mu_0(x)$. Introducing the following scaling:  $\alpha = \varepsilon, \eta=1/\varepsilon,$ for the binary interaction \eqref{eq:bin}, and  defining by $\mu^\varepsilon(x,t)$ a solution for the scaled equation  \eqref{eq:Boltz},
when $\varepsilon\to0$, $\mu^\varepsilon(x,t)$ converges point-wise, up to a subsequence, to  $\mu(x,t)$ satisfying  the following nonlinear mesoscopic equation
\begin{align}\label{eq:FP}
\pa_t \mu + \nabla_x\cdot\left((\mathcal{P}[\mu] + \bigK[\mu])\mu\right)= 0
\end{align}
with initial data $\mu_0(x)=\mu(x,0)$ and 
where
\begin{align}
&\mathcal{K}[\mu](x,t)  = \int_{\mathbb{R}^d}K(x,y)\mu(y,t)\,dy. \label{eq:kernelK}
\end{align} 
\end{thm}

We refer to \cite{ACFK16}, for a complete proof.
\begin{remark}
 Theorem \ref{thm:grazing} furnishes a consistency result to our approach, in particular in the quasi-invariant scale we can conclude that the  binary interaction dynamics converges to a sparse feedback constrained mean-field equation. 
 Moreover, substituting directly into equation \eqref{eq:FP} the empirical measures $\mu^N(t)$ concentrated onto $(x_i(t))_{i=1}^N$ for $t\geq 0$, restitutes exactly the MAS systems \eqref{eq:MAS2}, where the control $u_i(t)$ is expressed in feedback form as $u^*_i(t)= \sum_{j=1}^NK(x_i,x_j)/N$, for $i=1,\ldots,N$.
 \end{remark}

\section{Numerical approximation and tests}
\subsection{Numerical methods}
Standard schemes to solve Boltzmann-type equations are often based on the Monte Carlo method. For this, let us consider the initial value problem given by the equation \eqref{eq:Boltz}, in the quasi-invariant scaling, for which $\alpha =\epsi,\lambda = 1/\epsi$ and with initial data $\mu(x,0)=\mu_0(x)$, 
\begin{equation}\vspace{0.5em}\label{eq:Coll}
\dfrac{d}{dt}\mu(x,t) =  \dfrac{1}{\epsi}\left[{Q}_\epsi^{+}(\mu,\mu)(x,t)- \mu(x,t)\right],
\end{equation}
where $Q^{+}_\varepsilon(\cdot,\cdot)$ denotes the gain part, accounted for the the density of agents gained at state $x$ after the binary interaction  \eqref{eq:bin}, and the second term represents the loss. Next, we consider a first order forward scheme for the scaled Boltzmann-type equation \eqref{eq:Coll}
\begin{equation}\label{eq:MCBoltz}
        \mu^{k+1}=\left(1-\frac{\delta t}{\varepsilon}\right)\mu^{k}+\frac{\delta t}{\varepsilon}{{Q}_\varepsilon^{+}(\mu^k,\mu^k)},
\end{equation}
for $k=1,\ldots,M_{tot}-1$,  where we denote by $\mu^k$ the  approximation of $\mu(x,k\delta t)$, for $\delta t >0$ the time step discretization. Since $\mu^k$ is a probability density, thanks to mass conservation, also $Q_\varepsilon^{+}(\mu^k,\mu^k)$ is a probability density,  and therefore under the restriction $\delta t\leq\varepsilon$ then also $\mu^{k+1}$ is a probability density.
Note that, since we aim at small values of $\varepsilon$ and  under the condition $\delta t\leq\epsi$, the natural choice is to take $\delta t=\varepsilon$. This choice maximizes at every step the number of interactions among the agents. 

The approximation of the controlled system also requires a numerical realization of the optimal feedback controller. For the instantaneous case this follows from the closed form expression \eqref{IC}. For the infinite horizon sparse feedback control the numerical approximation of eq. \eqref{eq:DP} is performed by implementing a semi-Lagrangian, policy iteration scheme, following the design guidelines presented in \cite{AFK15,KKK16}.

Combining both methods, we  propose the following \emph{Binary Constrained Interaction} algorithm (BCI) to solve \eqref{eq:MCBoltz},
\begin{alg}[\textbf{BCI algorithm}]~\,
  \begin{enumerate}
  \item[\texttt{0.}] Compute the feedback control $S_\epsi(x,y)$ on a suitable discretization of the grid of $\Omega\times\Omega$.
  \item[\texttt{1.}]
  Given $N_s$ samples $\left\{x^0_k\right\}_{k=1}^{N_s}$, from the initial distribution $\mu_0(x)$;
  \item[\texttt{2.}] 
  \texttt{for} $k=0$ \texttt{to} $M_{tot}-1$
  \begin{enumerate}
  \item[\texttt{a.}]  set $N_c = \textsc{Iround}({N_s}/{2})$;
  \item[\texttt{b.}]  select $N_c$ random pairs $(i,j)$ uniformly without repetition among all possible pairs of individuals at time level $t_m$;
  \item[\texttt{c.}]  evaluate the interactions $P(x^k_i,x^k_j), P(x^k_j,x^k_i)$ and $S_\epsi(x^k_i,x^k_j), S_\epsi(x^k_j,x^k_i)$;
  \item[\texttt{d.}]  for each pair $(i,j)$, compute the post-interaction position $x_i^*$, $x_j^*$ via \eqref{eq:bin}.
  \item[\texttt{e.}]  set $x_i^{k+1}=x_i^{*}$, $x_j^{k+1}=x_j^{*}$.
  \end{enumerate}
\item[]\texttt{end for}
  \end{enumerate}
 \label{ANMC}
\end{alg}\medskip
Where the function $\textsc{Iround}(\cdot)$ denotes the integer stochastic rounding.
We refer to \cite{APb,PTa} for further details on this type of algorithms.

\subsection{Numerical tests}

We validate the consistency of our numerical procedure, considering the control of different social dynamics. We perform two tests: a first test for consensus dynamics  (\cite{hekr02}), and a second test for attraction-repulsion dynamics, see for example \cite{MR2257718}. In both cases $d=1$, and model parameters are reported in  Table \ref{tab:par} .
\begin{table}[h]
\caption{Test 1 and Test 2 parameters.}
\label{tab:all_parameters}
\begin{center}
\begin{tabular}{cccccc }
\hline
& $N_s$ & $\epsi $ & $T$  & $\bar{\gamma}$ & $\lambda$    \\
\hline
\hline
Test 1:  & $5\times10^5$ & $5\times10^{-5}$ & $40$ & $0.3$ & 0.05 \\
\hline
Test 2: & $5\times10^5$ & $5\times10^{-5}$ & $10$ &$0.25$ & 0.1  \\
\hline
\end{tabular}\label{tab:par}
\end{center}
\end{table}

Furthermore, in order to reconstruct the density of agents from the $N_s$ samples, we fix the finite sub-interval, $\Omega = [-1,+1]\subset\R$, with space  discretization step $\delta x  =  0.025$.
In both cases we compare the results of the sparse feedback with the optimal feedback controllers with $\ell_2$ control penalization,  as in \cite{AHP,ACFK16}, and we consider the control to be bounded in the set $U=[-1,1]$.

\subsubsection{Test 1: Hegselmann-Krause model.}
We consider the Hegselmann-Krause model, also known as bounded confidence model, which describes the evolution of opinion in a society. In particular,  the consensus process is weighted by the interaction function  $P(\cdot,\cdot)$ in \eqref{eq:MFC}, defined as 
\begin{align}\label{eq:MFSz2}
P(x,y)  = \psi^\epsilon_\Delta(|x-y|), 
\end{align}
where  $\psi^\epsilon_\Delta(r)$, for some $\epsilon> 0$ is a regularization of the characteristic function $\chi_{\{|x-y|\leq \Delta\}}(y)$. The parameter $\Delta>0$ represents the confidence level, namely the range in which the interaction among two agents can happen. 
This type of model describes the propensity of an agent with opinion $x$ to interact only within a confidence range $I_{\Delta} =[x-\Delta,x+\Delta]$. We fix $\Delta = 0.4$, and we study the evolution of the control problem up to time $T = 40$ with initial data uniformly distributed,  $\mu^0(x) \sim \textrm{Unif}([-1,1])$. In Figure \ref{Fig_T1}, we depict the evolution of the dynamics without any control, showing the emergence of two clusters.
\begin{figure}[t]
\centering
\includegraphics[scale=0.43]{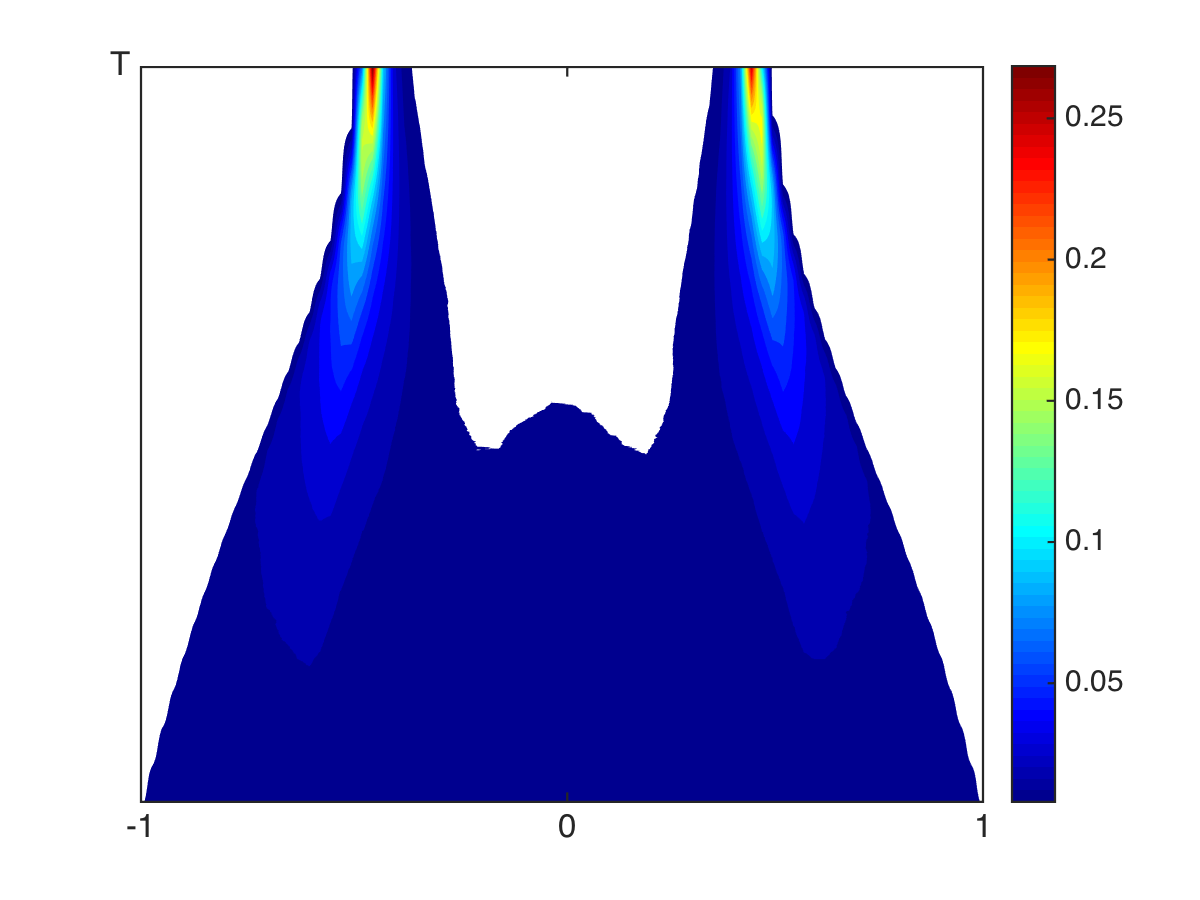}
\caption{Test \#1. Evolution of the Hegselmann-Krause model without control, with confidence level $\Delta = 0.4$.}\label{Fig_T1} 
\end{figure}

We introduce the action of feedback controls for the desired state $\hat x = 0$:  the instantaneous control (IC) in eq. \eqref{IC}, and the Hamilton-Jacobi feedback control (IH) in eq.  \eqref{eq:ochj}. We report in Figure \ref{Fig_T2}, the evolution of the constrained density, $\mu(x,t)$ and the evolution of the control action, $\mathcal{K}[\mu](x,t)$ in the instantaneous control case. The left-hand side depict the $\ell_2$ instantaneous control, whereas the right-hand side shows the $\ell_1$ instantaneous control. We remark the stronger action of the $\ell_2$ control with respect to a selective action of the control only outside a certain interval around the desired state $\hat x$. The evolution shows that in both cases consensus towards $\hat x$ is reached.
\begin{figure}[t]
\centering
\includegraphics[scale=0.215]{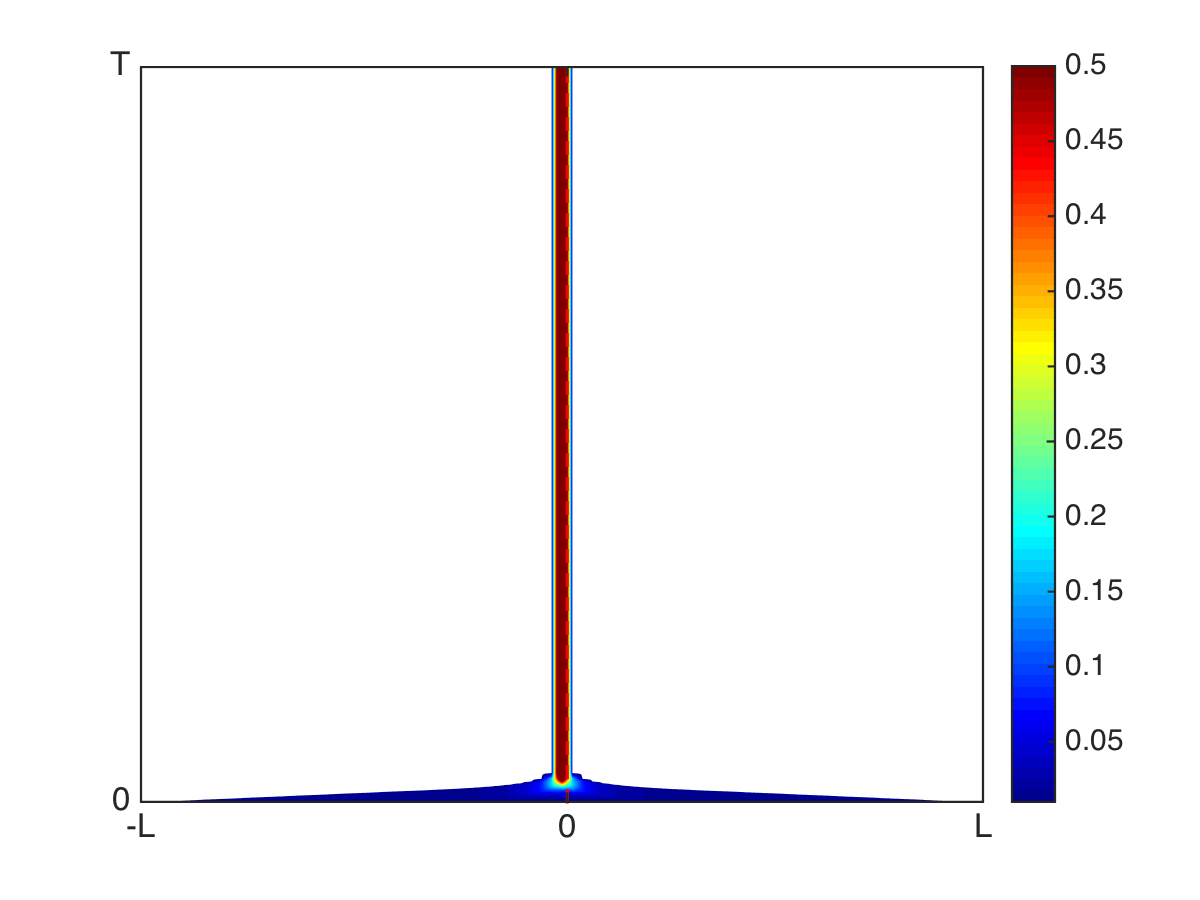}
\includegraphics[scale=0.215]{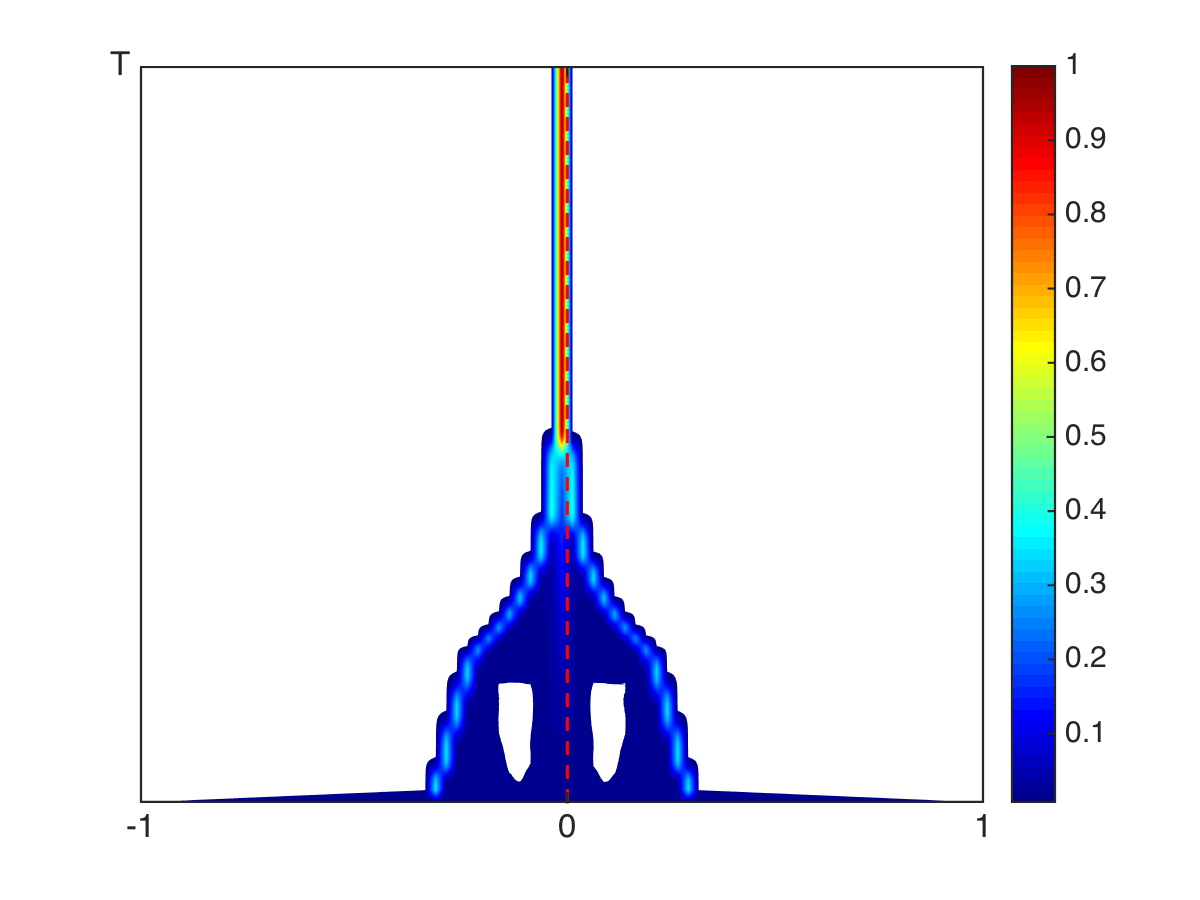}
\\
\includegraphics[scale=0.215]{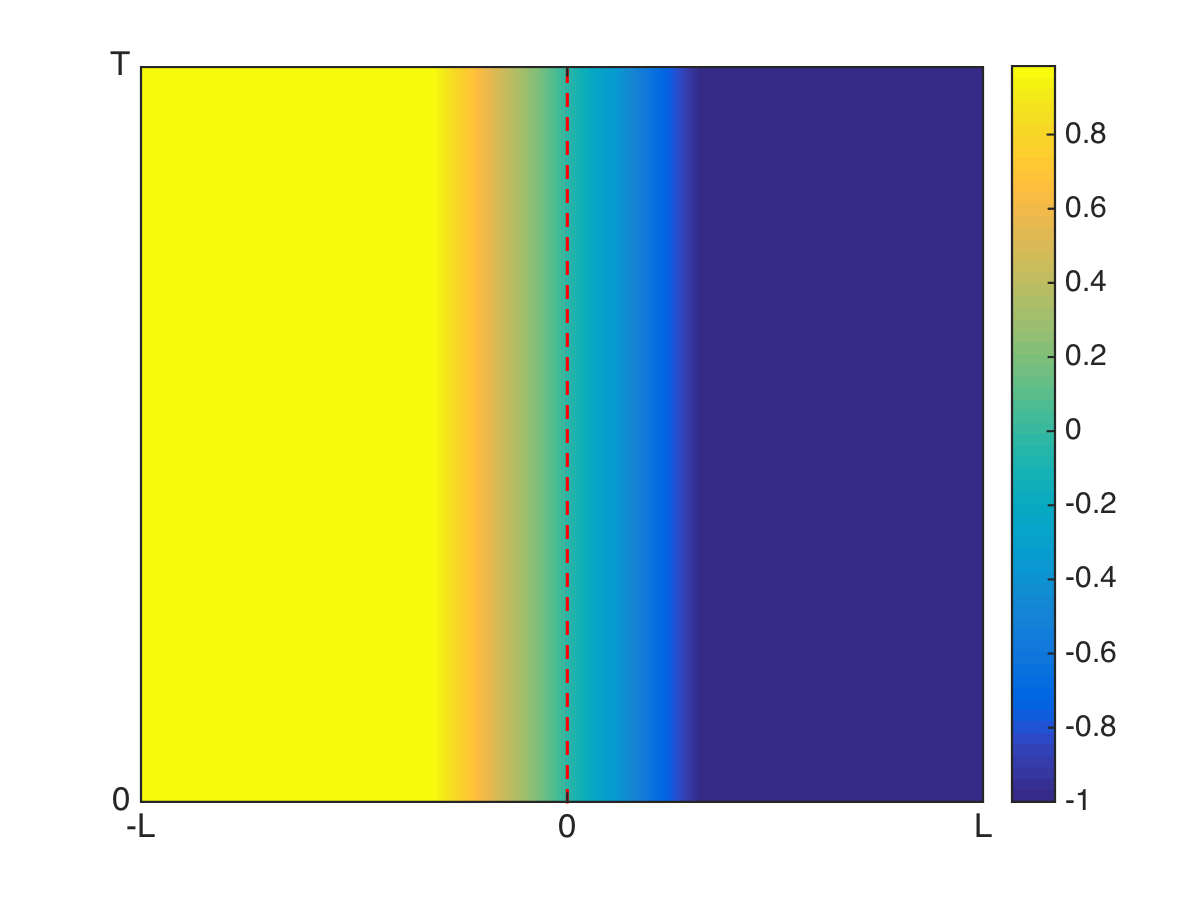}
\includegraphics[scale=0.215]{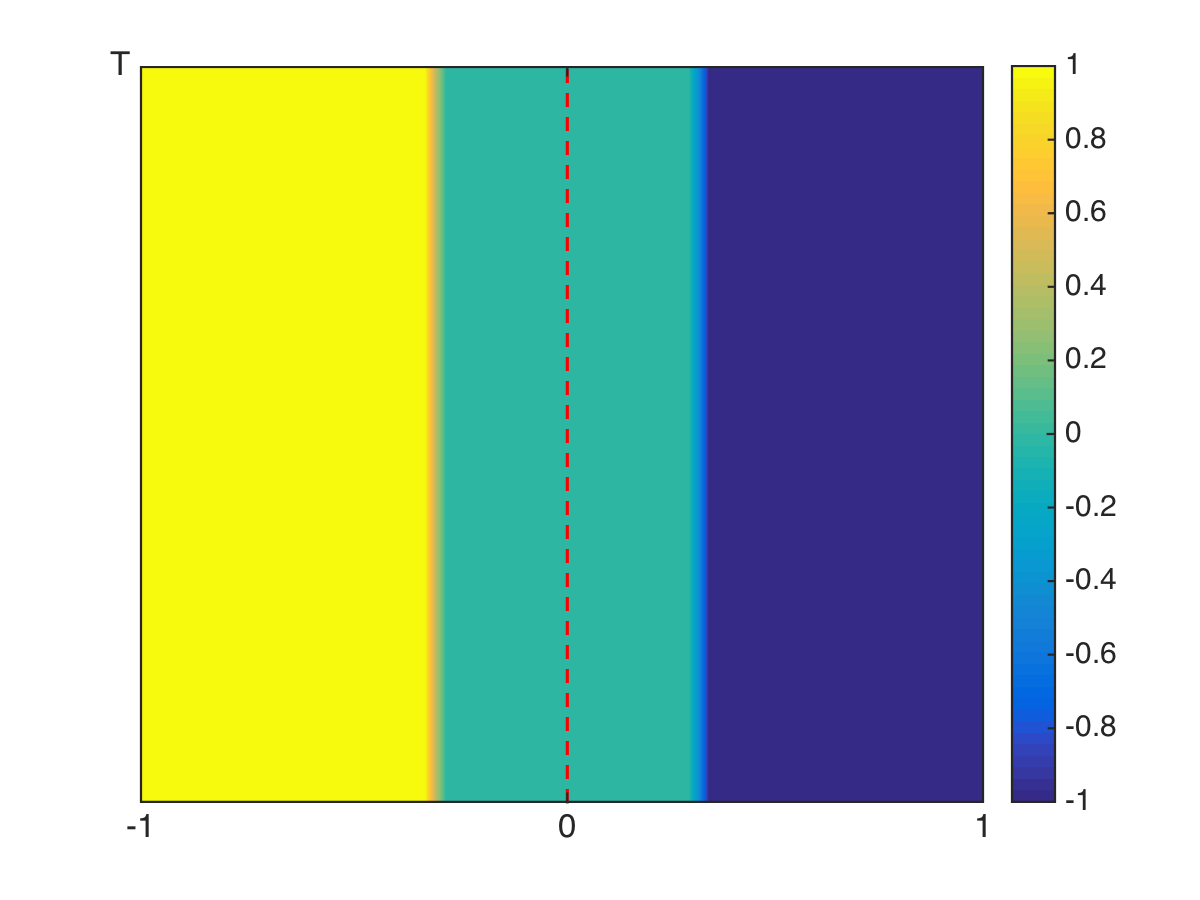}
\caption{Test \#1. Sparse instantaneous control (IC) in the  Hegselmann-Krause model. Left: $\ell_2$ feedback control. Right: $\ell_1$ sparse feedback control. Top: evolution of the density $\mu$ in the interval $[0,T]$.  Bottom: different actions of the control $\mathcal{K}[\mu]$, where the sparsity of the control is gained around the desired state $\hat x= 0$.}
\label{Fig_T2}
\end{figure}
Similarly Figure \ref{Fig_T3} depicts the action of the infinity horizon strategy (IH), in the $\ell_2$ and $\ell_1$ minimization setting, showing the convergence towards the desired state $\hat x$.

\begin{figure}[t]
\centering
\includegraphics[scale=0.215]{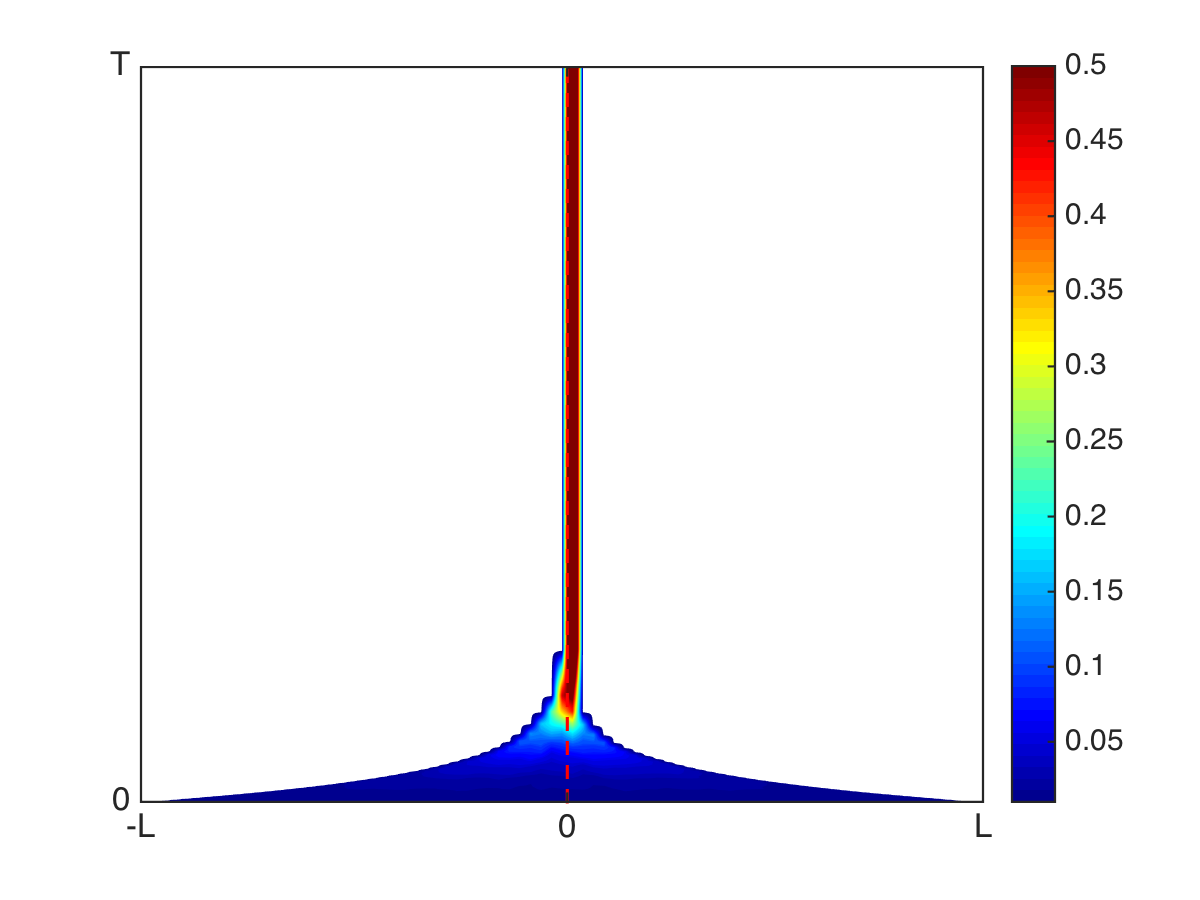}
\includegraphics[scale=0.215]{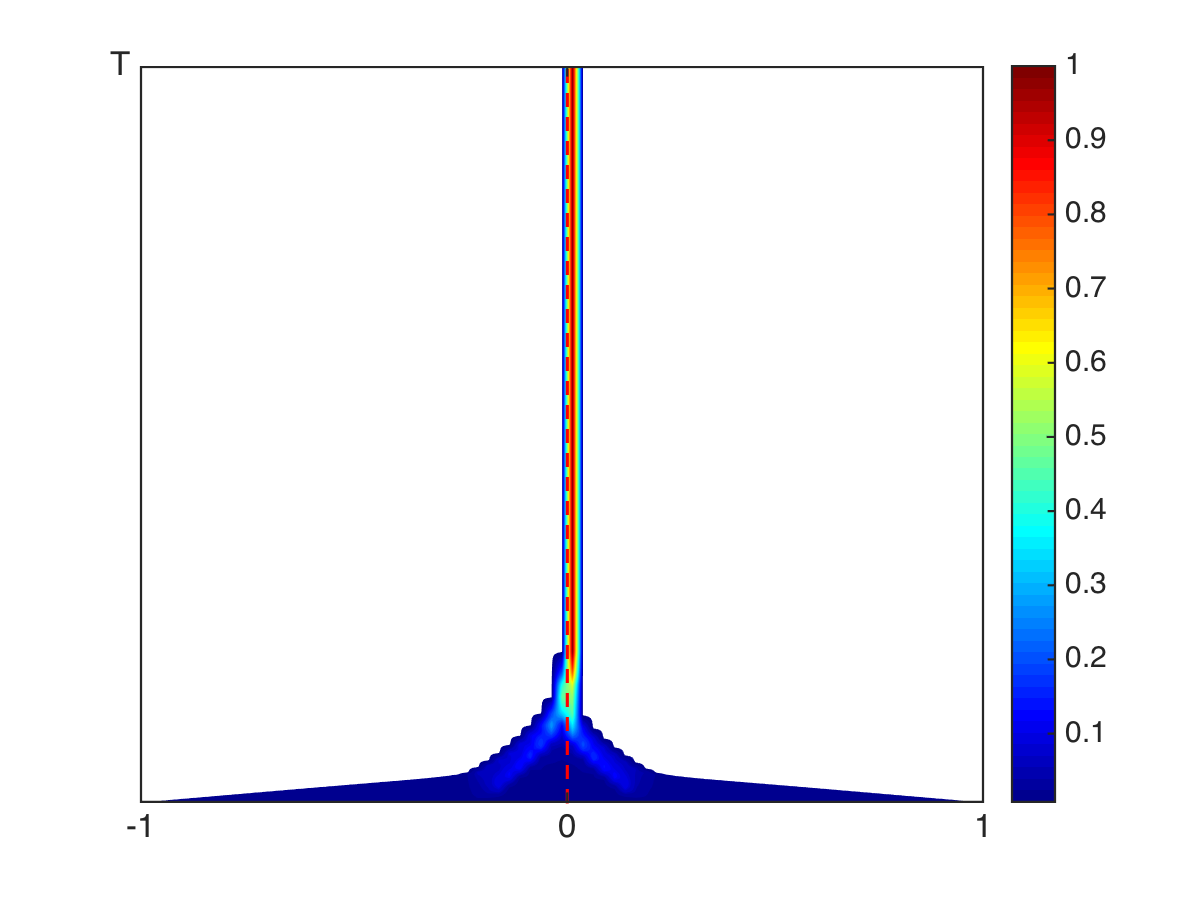}
\\
\includegraphics[scale=0.215]{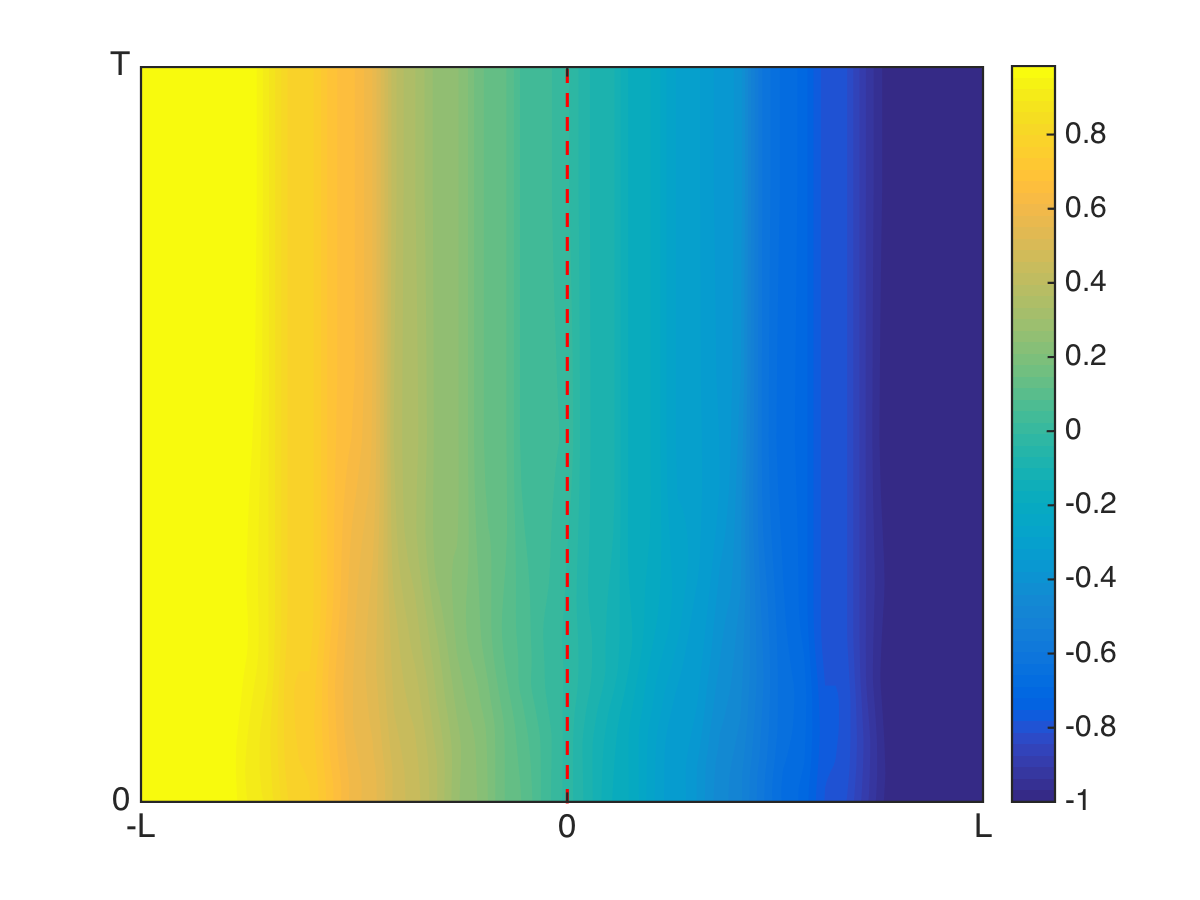}
\includegraphics[scale=0.215]{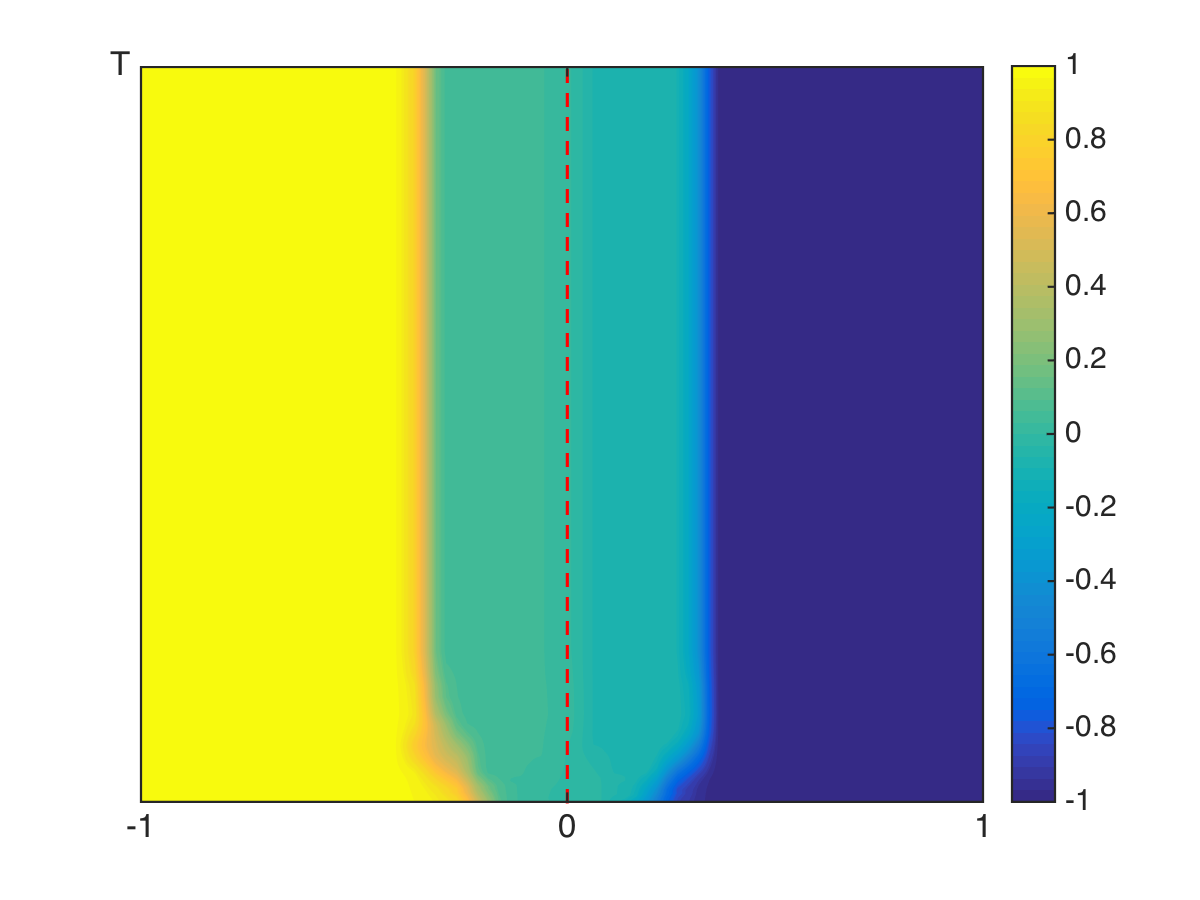}
\caption{Test \#1.  Infinite horizon control (IH) in the  attraction-repulsion model. On the left the $\ell_2$ feedback control, on the right the $\ell_1$ sparse control. Top:  evolution of the density $\mu$. Bottom:  evolution of the feedback control $\mathcal{K}[\mu]$.}
\label{Fig_T3}
\end{figure}

\subsubsection{Test 2: Attraction-Repulsion model.}
In this second test we consider a interaction potential $P(\cdot,\cdot)$ in \eqref{eq:MFC}, as a smoothed version of a power-law potential,
\begin{align}\label{eq:MFSz2}
P(x,y)  =(\sigma + |x-y|)^a - (\sigma+|x-y|)^b, 
\end{align}
where we fix $a = 1$, $b=-1$, respectively the power of the attractive part and the repulsive part, and $\sigma =10^{-4}$ is a regularization parameter.
We study the evolution of the control problem up to time $T = 10$ with initial data uniformly distributed,  $\mu^0(x) \sim \textrm{Unif}([-1,1])$. 
\noindent
We report in Figure \ref{Fig_AR0} the evolution of the uncontrolled dynamics, which shows a confinement of the density in a interval with higher concentration on the border of its support. Introducing the action of the control we want to steer the agents towards the $\hat x = 0 $ state. Similarly to the previous test, we compare the action of the $\ell_1$ sparse feedback control with the $\ell_2$-penalized control.
\begin{figure}[t]
\centering
\includegraphics[scale=0.4]{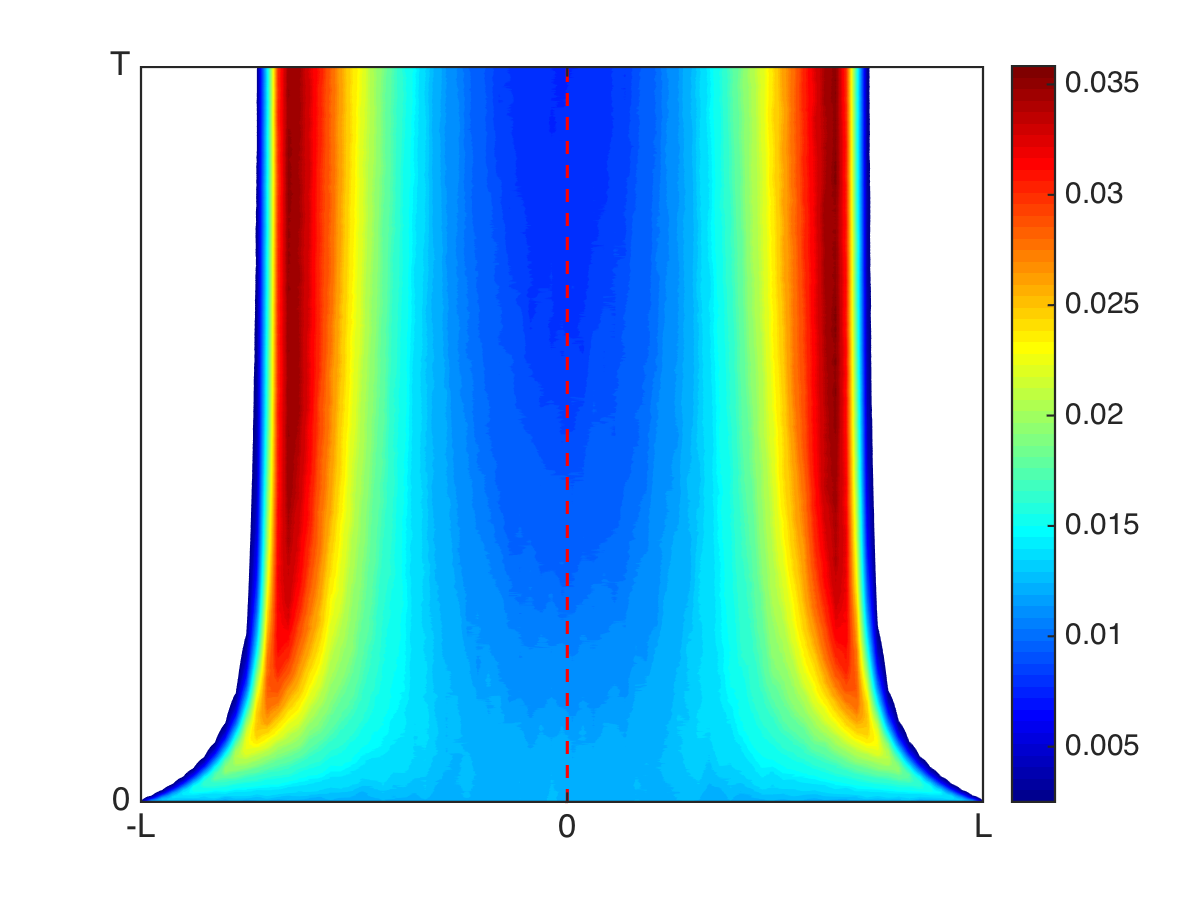}
\caption{Test \#2. Evolution of the attraction-repulsion model without control, with confidence level $\Delta = 0.4$.}
\label{Fig_AR0}
\end{figure}
In  Figure  \ref{Fig_AR1}, we depict the evolution of the dynamics for the instantaneous controls (IC), the test shows that, whereas the $\ell_2$ instantaneous control is capable to control the dynamics, the sparse instantaneous control fails to steer the system towards the reference state, and concentration of the density in two peaks appears. Alternatively, in Figure \ref{Fig_AR2}, we depict the evolution of the dynamics for the infinite horizon feedback control (IH), where the desired state is reached for both penalizations.
\begin{figure}[t]
\centering
\includegraphics[scale=0.215]{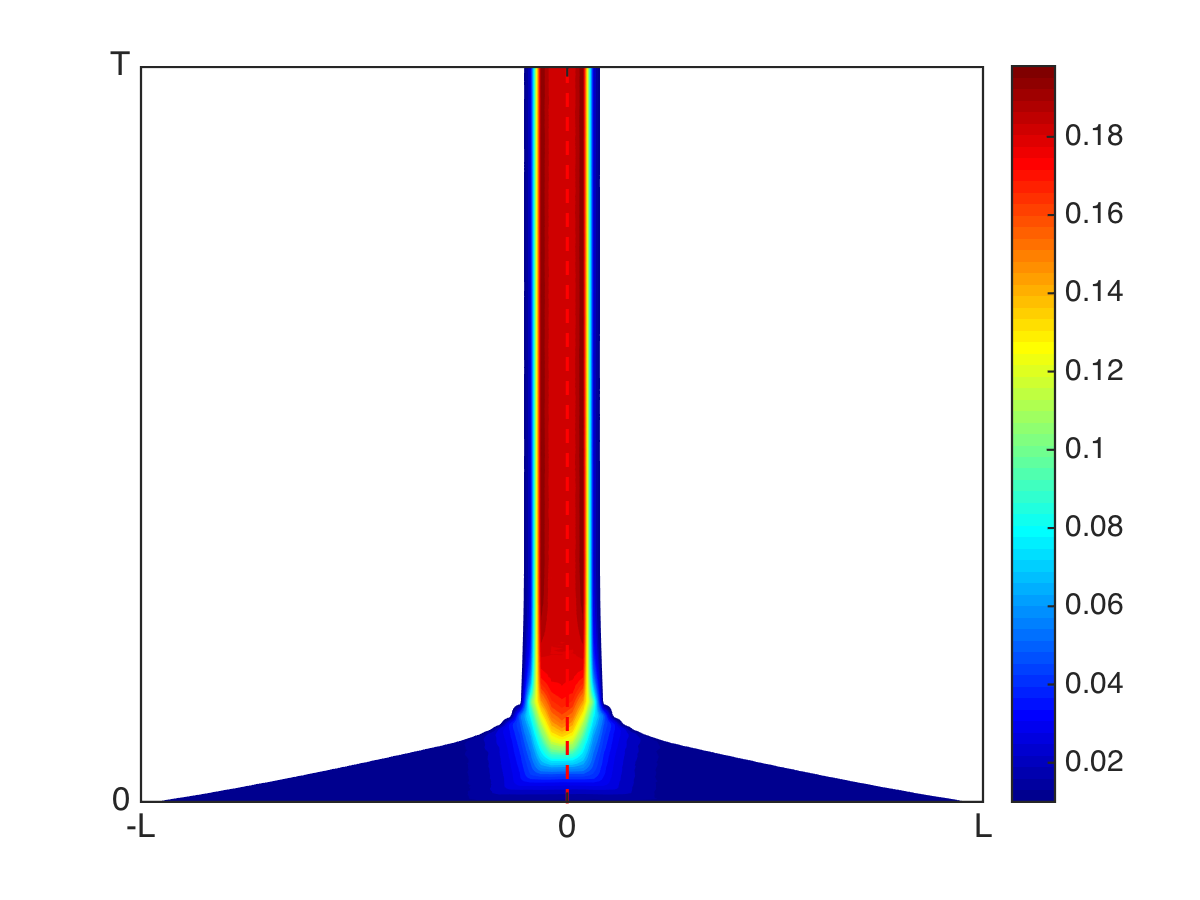}
\includegraphics[scale=0.215]{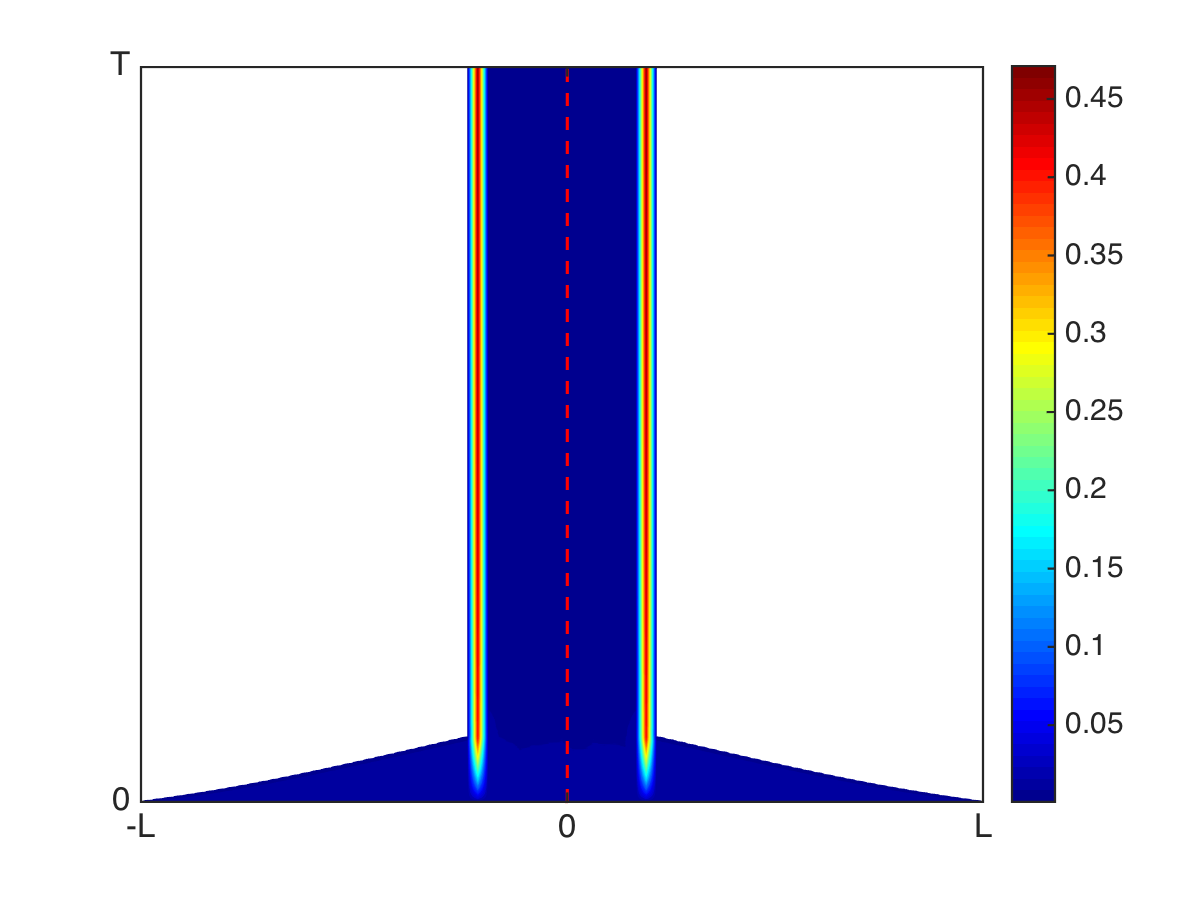}
\\
\includegraphics[scale=0.215]{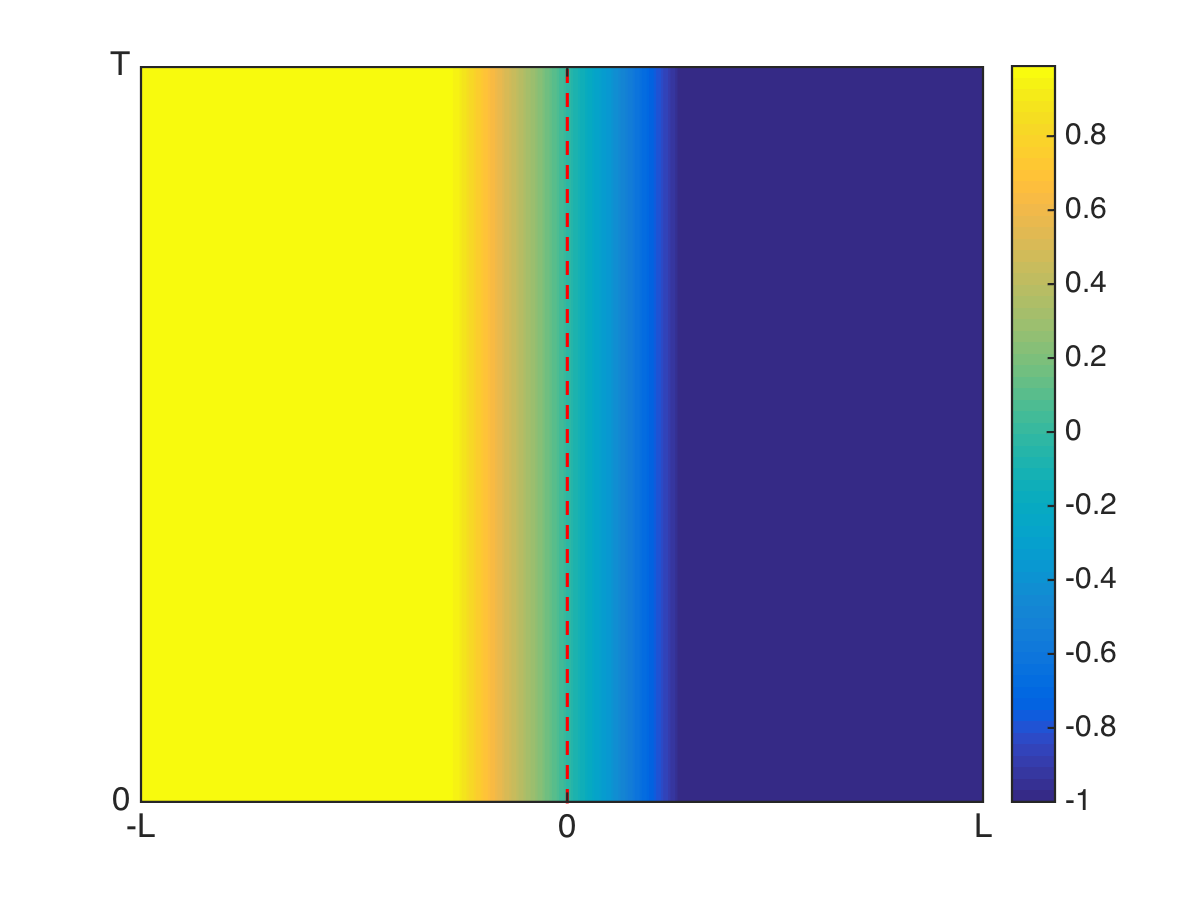}
\includegraphics[scale=0.215]{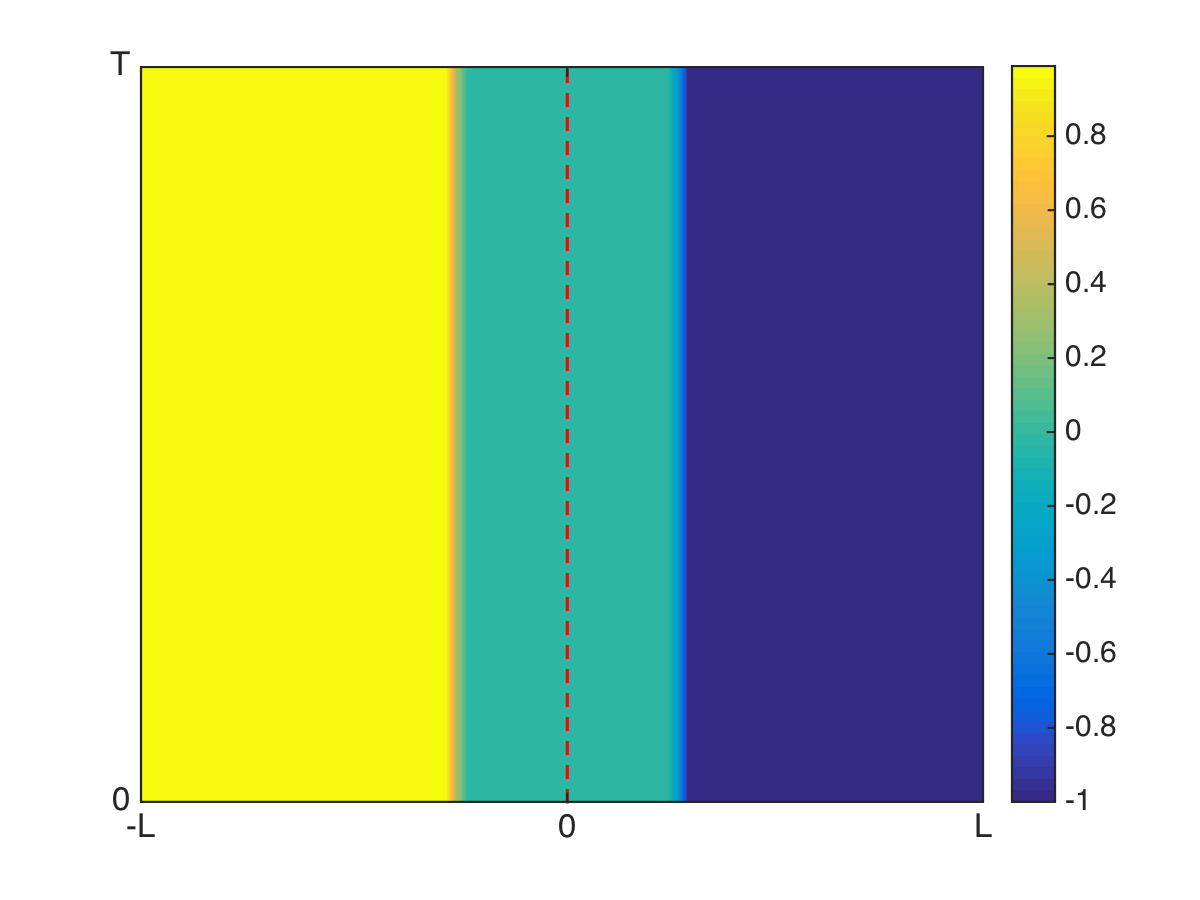}
\caption{Test \#2. Instantaneous control (IC) in the  attraction-repulsion model. On the left the $\ell_2$ minimization, on the right the $\ell_1$ minimization, first line reports the evolution of the density $\mu$, the second line the evolution of the feedback control $\mathcal{K}[\mu]$.}
\label{Fig_AR1}
\end{figure}

\begin{figure}[t]
\centering
\includegraphics[scale=0.215]{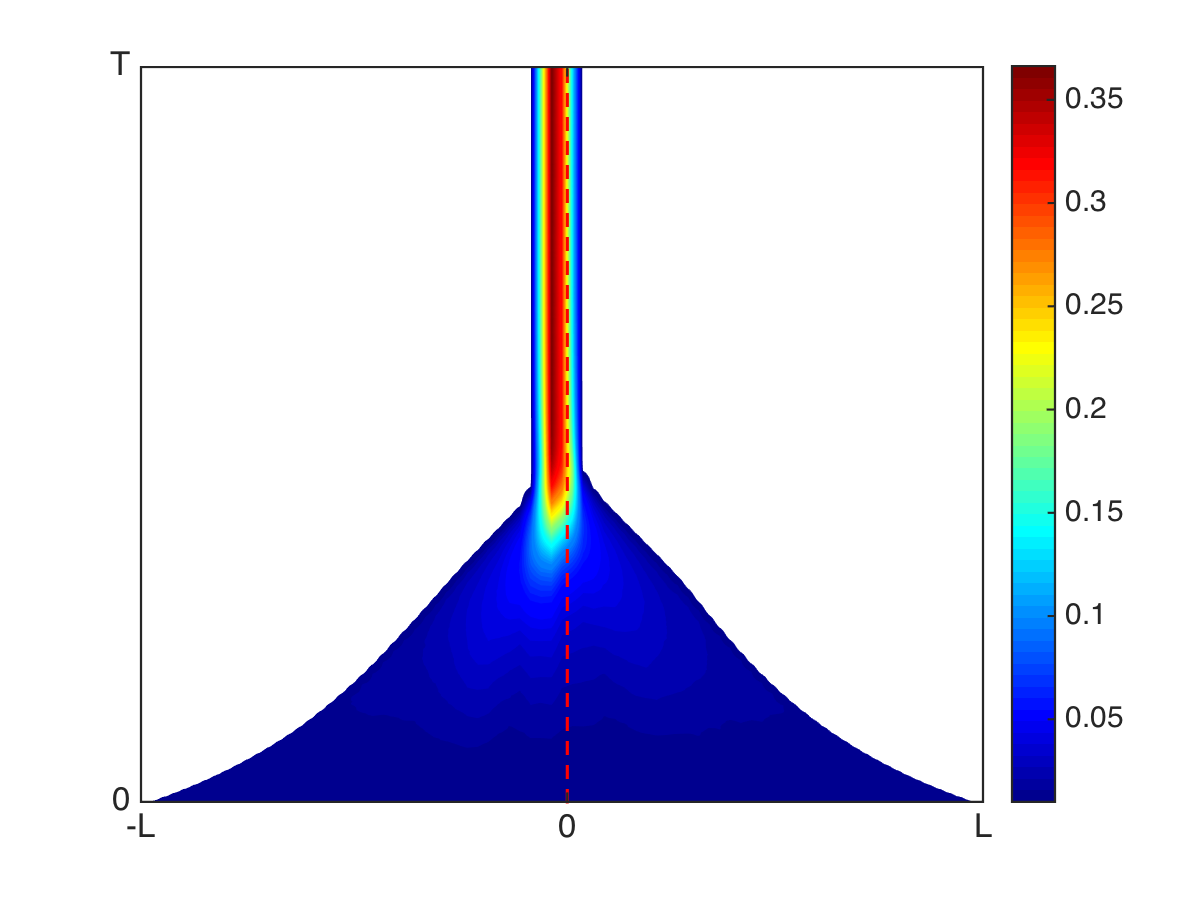}
\includegraphics[scale=0.215]{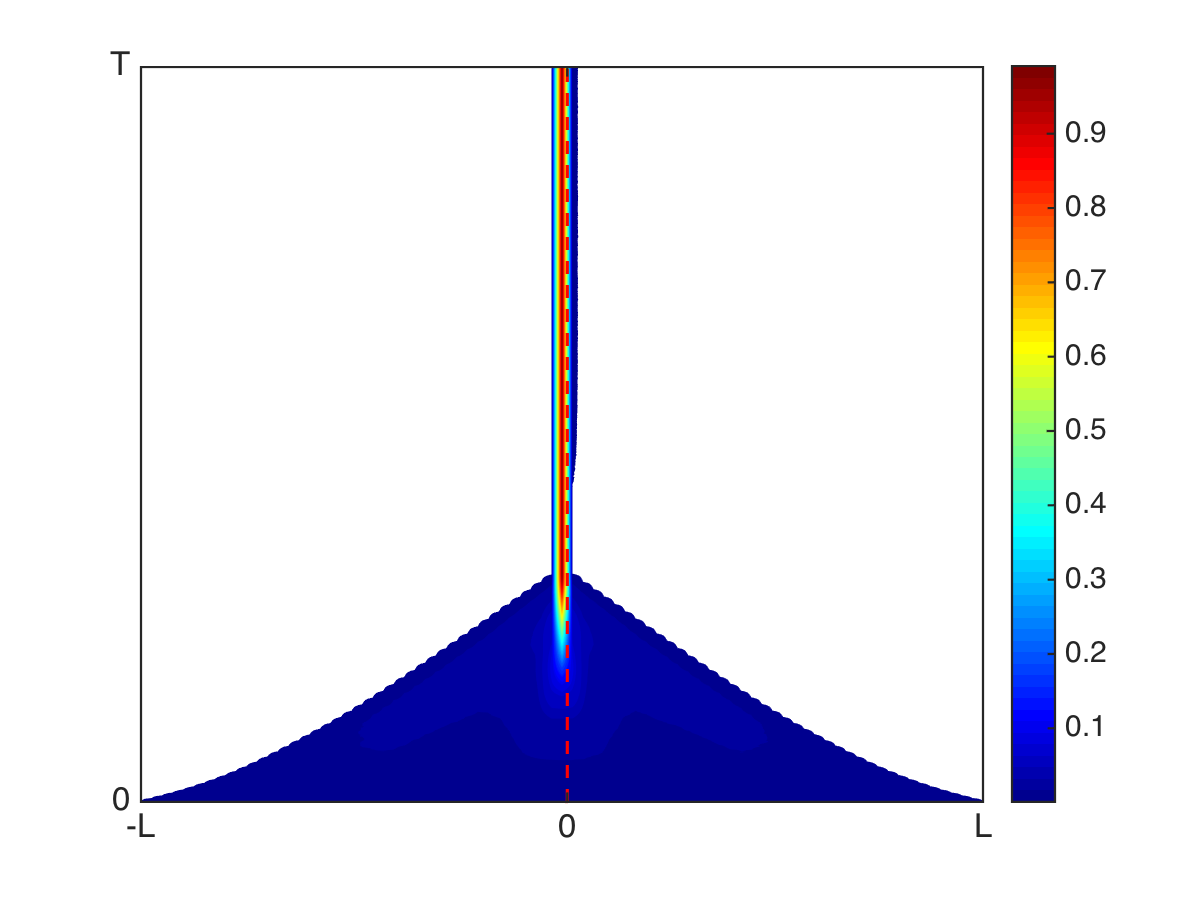}
\\
\includegraphics[scale=0.215]{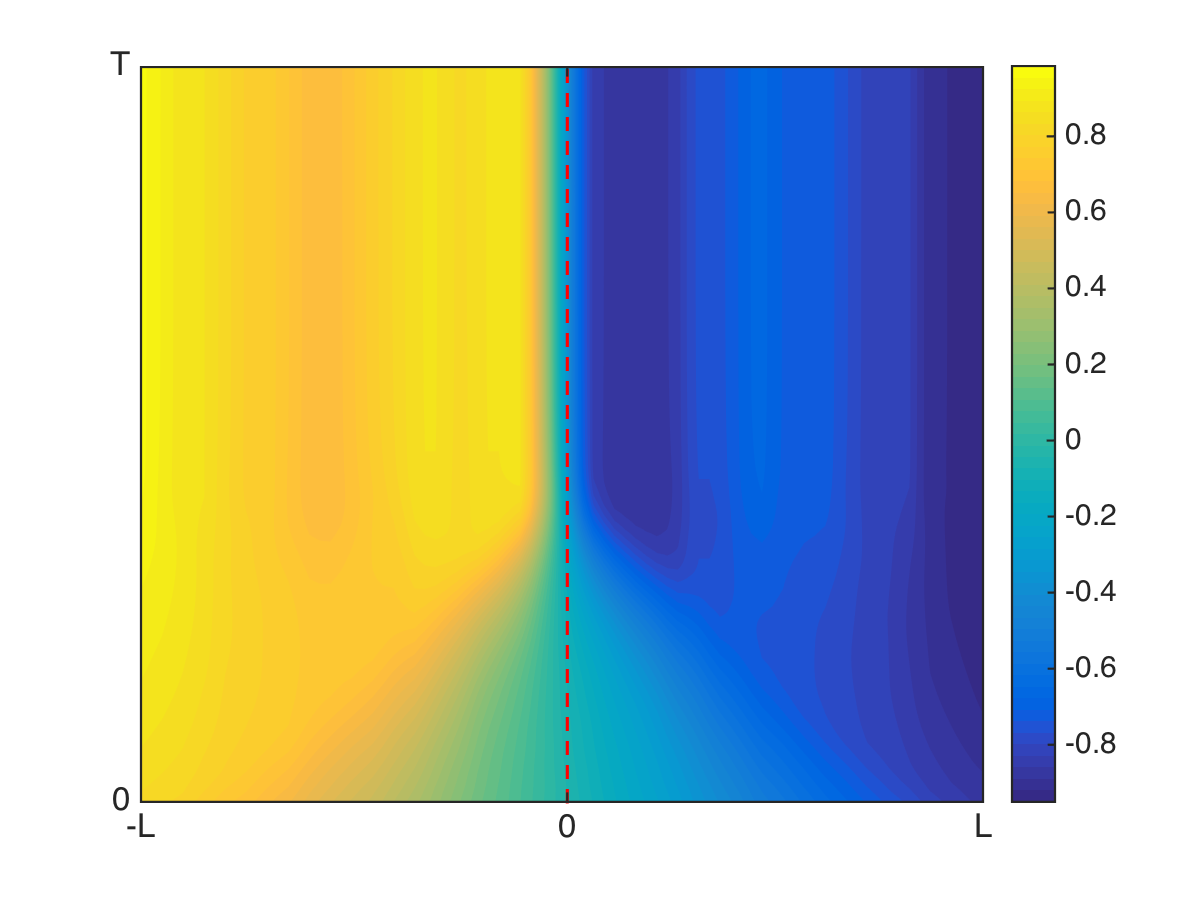}
\includegraphics[scale=0.215]{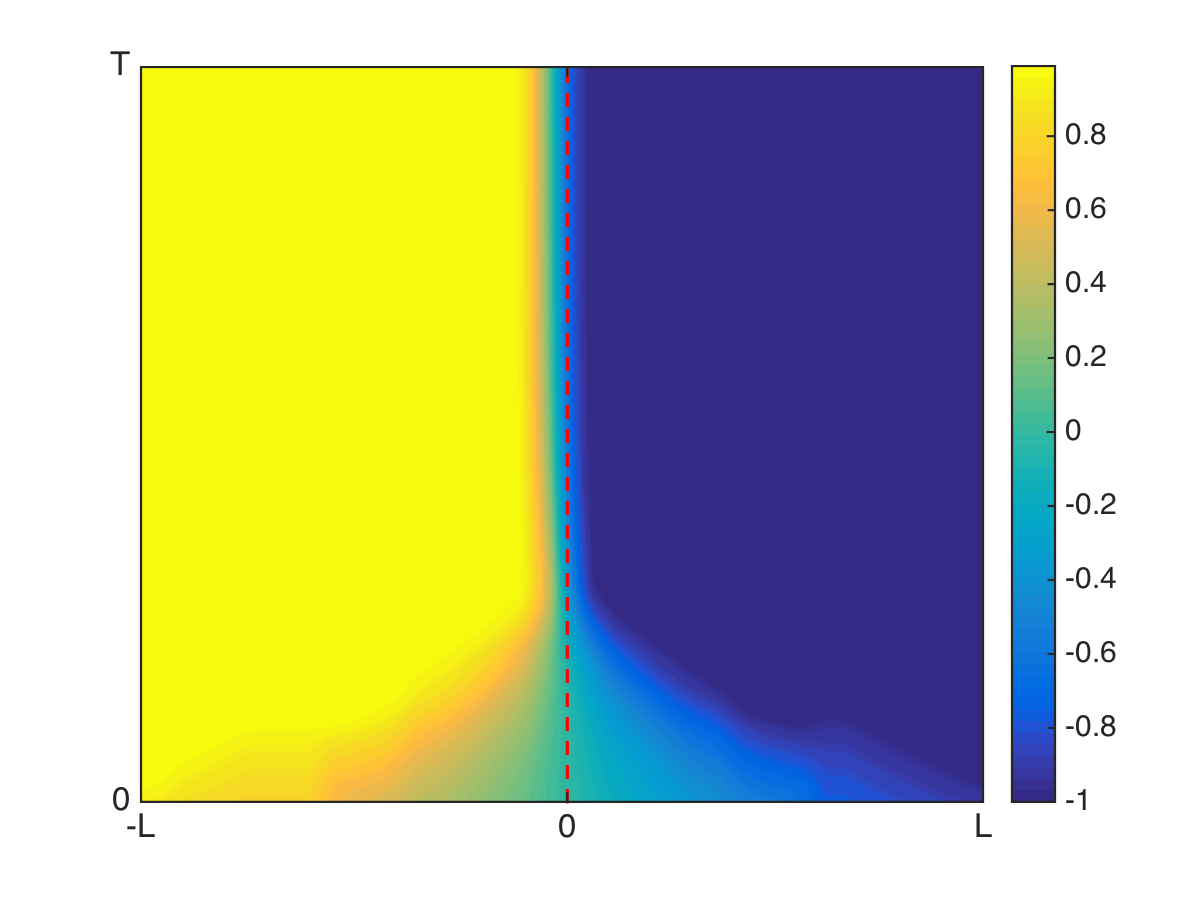}
\caption{Test \#2. Infinite horizon control (IH) in the  attraction-repulsion model. On the left $\ell_2$ minimization, on the right the $\ell_1$ minimization, first line reports the evolution of the density $\mu$, the second line the evolution of the feedback control $\mathcal{K}[\mu]$.}
\label{Fig_AR2}
\end{figure}


\end{document}